\input amstex
\input amsppt.sty
\input epsf
\magnification\magstep1
\pageheight{7.5in}
\pagewidth{5.5in}
\NoBlackBoxes

\def\R{\text{\bf R}}

\topmatter
\author
P. Duvall, J. Keesling 
\endauthor

\title
The Hausdorff dimension of the boundary of the L\'evy dragon
\endtitle

\abstract 
A theoretical approach to computing the Hausdorff dimension of the 
topological boundary of attractors of iterated function systems is
developed. The curve known as the L\'evy Dragon is then studied in
detail and the Hausdorff dimension of its boundary is computed 
using the theory developed. The actual computation is a complicated
procedure.  It involves a great deal of combinatorial topology as 
well as determining the structure and certain eigenvalues of a
$752 \times 752$ matrix.  Perron-Frobenius theory plays an important role
in analyzing this matrix.
\endabstract

\thanks This research was done while the first author was a Visiting
Professor at the University of Florida. He gratefully acknowledges UFL's 
hospitality, as well as research
support from the University of
North Carolina, Greensboro
\endthanks

\address 
Department of Mathematical Sciences,
University of North Carolina at Greensboro,
Greensboro, NC 27412, USA
\endaddress
\email duvallp\@uncg.edu
\endemail
\address 
University of Florida, Department of Mathematics, P.O.~Box~118105,
358 Little Hall, Gainesville, FL 32611-8105, USA
\endaddress
\email  jek\@math.ufl.edu
\endemail

\subjclass Primary 28A20, 54E40, 57Nxx
\endsubjclass

\keywords  Hausdorff dimension,
iterated function systems,
attractors, fractal geometry
\endkeywords

\endtopmatter

\document

\heading
\S1. Introduction
\endheading

In \cite{L1}, P. L\'evy studied the complex curve which has come to be known as 
the L\'evy Dragon. He showed that the Dragon is self-similar, and argued that
the plane can be tiled by copies of it. We give a self-contained account of
the Dragon that determines the properties of this set that we need in our
computations.  However, the paper by L\'evy is a {\it tour de force} on this
subject.  It is amazing that he was able to determine so many properties of
this and other curves in \cite{L1} without the use of modern computers.
The reader is urged to see \cite{L2} for an 
excellent annotated translation of \cite{L1}.  
\medskip

\epsfysize=1.8in
\centerline{\epsffile{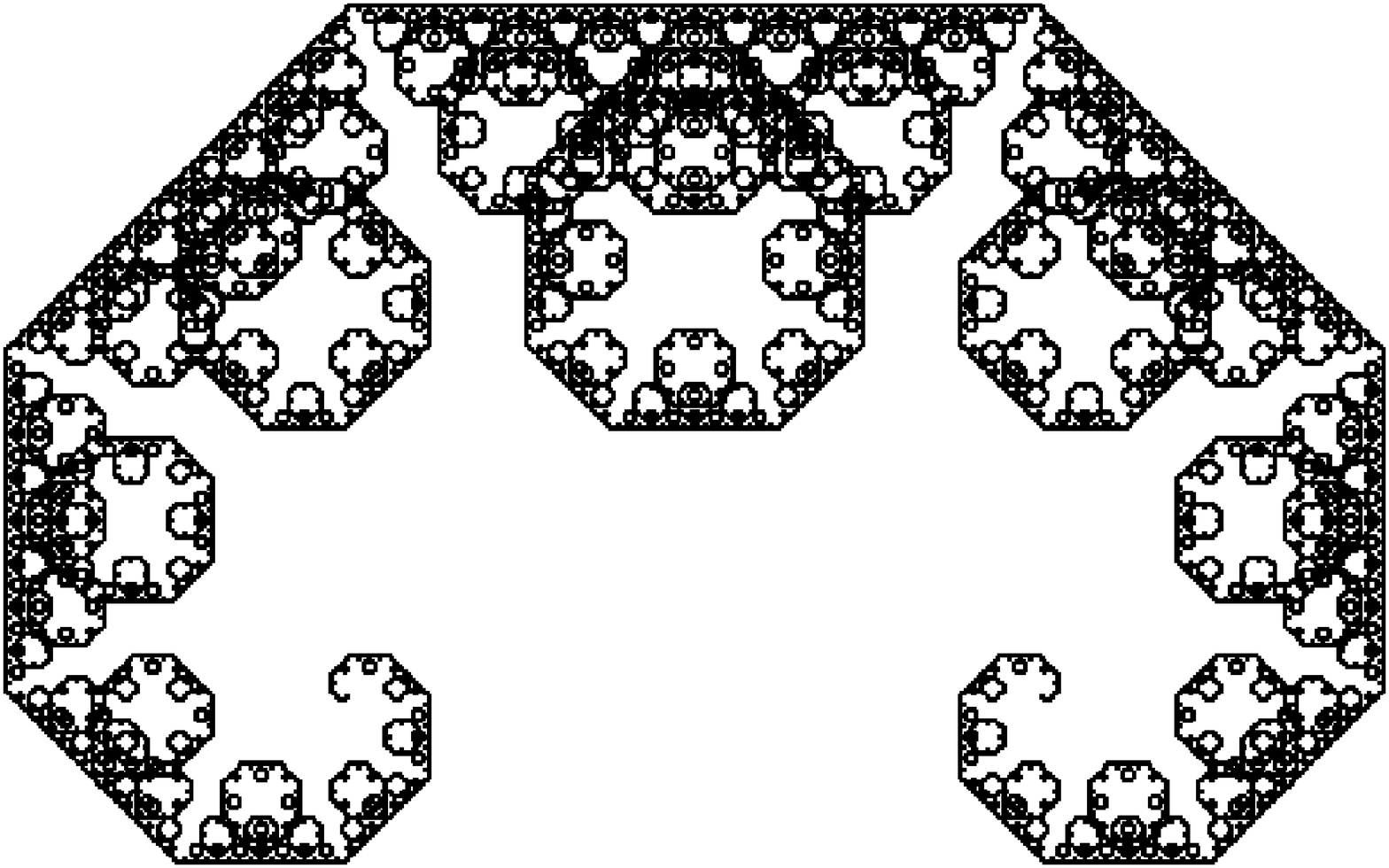}}
\smallskip
\centerline{Figure 1. The L\'evy Dragon}
\medskip

        Although L\'evy found so many properties of the Levy Dragon,
he did not compute the Hausdorff dimension of its boundary.  In
\cite{L2, page 236} G. Edgar asks whether the dimension of the boundary
of the Dragon is greater than one.  In this paper we answer that
question in the affirmative and in fact determine a precise value
for this dimension.  In the process of doing this we develop a
general theory for computing the Hausdorff dimension of the
boundary of self-similar fractal tiles.
This depends heavily on the paper by K. Falconer \cite{F1} who
gave an (abstract) method for determining the Hausdorff dimension of 
sub-self-similar sets.  Our contribution is to show that in the case of
boundaries, it is possible to derive a method which can be carried to a 
computation.  We then make a detailed study of the L\'evy Dragon
and show how to apply our theory to compute the dimension of its
boundary.  In subsequent work the authors together with A. Vince
have found a general method for finding the dimension of the boundary
of self-similar digit tiles \cite{DKV}.  The size of the matrices 
involved in computing
the dimensions of the boundaries of these tiles is typically less than 
$10 \times 10$.
Unfortunately the
L\'evy Dragon is not a self-similar digit tile and so these methods
do not apply directly to this example.  Fortunately, motivated by some
ideas in this paper R. Strichartz and Y. Wang \cite{SW} found a method
of computing the dimension of the boundary of the L\'evy Dragon
which has considerably reduced the computations required in this case as well. 
In their method
one only needs to determine a certain eigenvalue of an $11 \times 11$ 
matrix.  Strichartz and Wang made use of the graph directed constructions
of D. Mauldin and S. Williams \cite{MW} rather than the sub-self-similar sets
of Falconer \cite{F1} used in this paper and in \cite{DKV}.

First, we need to establish some terminology and notation. We assume that the
reader is familiar with the basic theory of self-similar fractals.  The books
by Edgar \cite{E1} and Falconer \cite{F2} are excellent references. The book
\cite{E2} contains a number of papers which give a nice historical 
perspective on
the subject.

If $X$ is a complete metric space and $F = \{f_{1}, f_{2}, \cdots,
f_{n}\}$ is a collection of contraction mappings of $X$ to itself,
then $F$ is said to be a (hyperbolic) {\it Iterated Function
System} or {\it IFS}. It is well known \cite{H}  that for  such an $F$
there is a unique compact set $K$ such that $K = \bigcup_{i=1}^{n}
f_{i}(K).$ $K$ is called the {\it attractor\/} or 
{\it invariant set} for $F$. In fact, if ${\Cal C}(X)$
denotes the collection of non-empty compact subsets of $X$
with the Hausdorff metric, then the map
$F:{\Cal C}(X) \rightarrow {\Cal C}(X)$ given by 
$F(A)=\bigcup_{i=1}^{n}f_{i}(A)$ is a
contraction mapping.  Thus $K$ is the
unique fixed point of $F$ guaranteed by Banach's contraction mapping
principle, and ${\displaystyle K=\lim_{k \rightarrow \infty} F^k(A)}$
for any $A \in {\Cal C}(X)$.  In this paper, we
will only be concerned with the case in which the $f_i$ are
{\it contracting similitudes\/}, maps satisfying 
$d(f_i(x),f_i(y))=c_id(x,y)$, where $d$ is the
metric on $X$ and $0<c_i<1$. 

Our main interest in this paper is the L\'evy Dragon.
It is the attractor for the
$IFS$ $F=\{f_1,f_2\}$, where $f_1$ and $f_2$ are similitudes of Euclidean space
$\R^2$ given by 
$$\align
 f_1(x,y)&= (\frac{x-y}{2},\frac{x+y}{2}), \\
 f_2(x,y)&= (\frac{x+y+1}{2},\frac{y-x+1}{2}).
\endalign$$
However, much of the theory developed in the early sections holds in 
a much broader setting.  So,
it is appropriate to have fairly general notational conventions.
The Hausdorff dimension of a set $E$ is denoted by $\dim_H E$, and the lower
and upper box-counting dimensions of $E$ are denoted by 
$\underline{\dim}_B E$ and $\overline{\dim}_B E$. Definitions of of these 
dimension functions can be found in \cite{F1}, \cite{F2} and \cite{E}.

Let $\Omega$ denote the collection of sequences
$I=\{i_1,i_2, \cdots\}$ with $i_k \in \{1,2,\cdots,n\}$. 
$\Omega_k$ will denote the set of sequences of length $k$ with 
entries from $\{ 1,2,\cdots ,n\}$.
 $\Omega$ has the natural metric given
by $$d(I,J)=\cases 1,&\text{if $i_1 \neq j_1$}\\
  c_{i_1}c_{i_2}\cdots c_{i_m},&\text{if $i_k=j_k$ for $k \leq m$ and 
$i_{m+1} \neq
j_{m+1}$.}\endcases$$ For $I \in \Omega$, we let $I_k$ denote the finite
sequence $\{i_1,i_2,\cdots,i_k\}$, and will often use the shorthand 
$f_{I_k}$ for the composition $f_{i_1}f_{i_2}\cdots f_{i_k}$ and $c_{I_k}$
for the product $c_{i_1}c_{i_2}\cdots c_{i_k}$. We will also need the
maps $\sigma :\Omega \rightarrow \Omega$ and $g:\Omega \rightarrow K$ given
by $\sigma ({i_1,i_2,\cdots})=\{i_2,i_3,\cdots\}$ and 
${\displaystyle g(I)=
\lim_{k \rightarrow \infty}f_{i_1}f_{i_2}\cdots f_{i_k}(x)}$. It is not
difficult to see that $g(I)$ is independent of the choice of $x \in X$.  
The map $\sigma$ is the well-known shift map on the sequence space.

\heading
\S2. Sub-Self-Similar Sets
\endheading

In \cite{F1}, Falconer introduced the notion of sub-self-similar (s.s.s.) 
sets,
and gave an (abstract) way to determine the dimension of such sets. 
In this section, we give a brief summary, in our notation, of the notions
from \cite{F1} that we will need for our calculations.
Given a
set $F=\{f_{1}, f_{2}, \cdots,f_{n}\}$ of contracting similitudes on $\R^l$,
the closed set $E$ is said to be {\it sub-self-similar} for $F$ if 
$E \subset \bigcup_{i=1}^{n}f_{i}(E)$. It follows that $E \subset K$, where
$K$ is the attractor for $F$. Note in particular that the topological 
boundary $\partial K$ of $K$ is s.s.s., since the interior int $K$ of $K$ is
mapped into itself by each of the open mappings $f_i$. The first crucial 
observation about s.s.s. sets is 
\proclaim{2.1 Proposition} Let $E$ be a closed set.  Then $E$ is s.s.s. for 
$F$ if and only if $E=g(A)$ for
some compact set $A \subset \Omega$ such that $\sigma (A) \subset A$.
Furthermore, if $E$ is s.s.s., such a set $A$ is given by 
$A=\{I \in \Omega| g(\sigma^k(I)) \in E$ for all $k \geq 0\}$.
\endproclaim
\demo{Proof}
This is Proposition 2.1 (and its proof) in \cite{F1}.\qed
\enddemo 
The set $A$ in Proposition 2.1 is fundamental in the calculation of 
$\dim_H E$. For $k>0$, let $A_k$ be the set of finite sequences obtained
by truncating elements of $A$ after $k$ terms. For $s \geq 0$ define
$$ \tau (s) = \lim_{k \rightarrow \infty} \left( \sum_{I_k \in A_k}
c_{I_k}^s \right)^{1/k}.$$ Recall that $F$ satisfies the {\it open
set condition} \cite{H} if there is a bounded nonempty open set 
$U$ such that 
$$ \align
\bigcup_{i=1}^n f_i(U) \subset U &, \text{and}\\
f_i(U) \cap f_j(U) = \emptyset &, \text{if } i \neq j.
\endalign$$
Our main tool for calculating dimensions is the following theorem of 
Falconer.

\proclaim{2.2 Theorem} Let $E$ be a closed set which is 
s.s.s. with respect to the family 
$\{f_{1}, f_{2}, \cdots,f_{n}\}$ of contracting similitudes which satisfies
the open set condition.  Let $s$ be the number satisfying $\tau (s) = 1$.
Then $s = \dim_H E = \underline{\dim}_B E = \overline{\dim}_B E$.
\endproclaim
\demo{Proof}
This is Theorem 3.5 of \cite{F1}.
The packing dimension $dim_P(X)$ is currently an important dimension function as well.
Since $dim_H(E) \leq dim_P(E) \leq dim_B(E)$ is always satisfied, we also have that
$s = dim_P(E)$ as well.  \qed
\enddemo

\heading
\S3. Boundaries of Attractors
\endheading

In this section, we specialize to the case of the boundary of the attractor
of an $IFS$ satisfying the open set condition, and show that stronger
invariance properties of the boundary give a more tractable version of
Theorem 2.2 for computing its Hausdorff dimension. Let 
$F=\{f_{1}, f_{2}, \cdots,f_{n}\}$ be an $IFS$ satisfying 
the open set condition, let $K$ be
its attractor, and let $A \subset \Omega$ be as in the previous section. One
obvious difficulty in applying Theorem 2.2 is that one has to sum over finite
sequences determined by an infinite condition. 
Let $$A'= \{I \in \Omega| f_{I_n}(K) \cap \partial K \neq \emptyset
  \text{ for all } n > 0\}.$$ We will show that $A'=A$ and use that
 to show that
membership in $A_n$ can be determined from finite data.

\proclaim{3.1 Proposition} $\sigma (A') \subset A'$.\endproclaim
\demo{Proof} Suppose that $I \subset A'$ but $\sigma(I) \notin A'$. Then for
some $n$, $$f_{i_2}f_{i_3}\cdots f_{i_n}(K) \subset \text{int }K.$$ But since
each $f_i$ maps the interior of $K$ into itself, 
$$f_{i_1}f_{i_2}\cdots f_{i_n}(K) \subset f_{i_1}(\text{int }K)
\subset \text{int }K,$$ a contradiction. \qed
\enddemo

\proclaim{3.2 Proposition} $A = A'$. \endproclaim
\demo{Proof}
Suppose that $I \in A'$. Let $k$ be given, and use $p \in K$ to compute $g$.
For any $m>0$, 3.1 implies that $f_{i_k}f_{i_{k+1}}\cdots f_{i_{k+m}}(K)
\cap \partial K \neq \emptyset$, 
so $$d(f_{i_k}f_{i_{k+1}}\cdots f_{i_{k+m}}(p),
\partial K) \leq c_{i_k}c_{i_{k+1}}\cdots c_{i_{k+m}}.$$ Therefore
${\displaystyle \lim_{m \rightarrow \infty} f_{i_k}f_{i_{k+1}}\cdots 
f{i_{k+m}}(p) \in \partial K}$ and $I \in A$ and $A' \subset A$.

Suppose that $I \in A$, but $f_{i_1}f_{i_2}\cdots f_{i_k}(K) \cap 
\partial K = \emptyset$. Then $$f_{i_1}f_{i_2}\cdots f_{i_k}(K) \subset 
\text{int }K.$$ Since 
$$f_{i_1}f_{i_2}\cdots f_{i_k}(K) \supset f_{i_1}f_{i_2}\cdots f_{i_{k+1}}(K)
\supset \cdots \supset f_{i_1}f_{i_2}\cdots f_{i_{k+m}}(K)$$ for all $m$, 
it follows that $g(I) \in \text{int }K$, a contradiction. 
Thus $A \subset A'$. \qed
\enddemo

\proclaim{3.3 Proposition} $A_k=\{I \in \Omega_k|f_{I_k}(K) \cap \partial K
\neq \emptyset \}$.
\endproclaim
\demo{Proof}
Let $I_k$ be a sequence in $\Omega_k$ such that $f_{I_k}(K) \cap \partial K
\neq \emptyset$.
If $f_{I_{k-1}}(K) \subset \text{int } K$, we would have that 
$f_{I_k}(K) = f_{I_{k-1}}(f_{i_k}(K)) \subset f_{I_{k-1}}(K) 
\subset \text{int } K$, a contradiction. Thus $f_{I_{k-1}}(K) \cap K
\neq \emptyset$, and, working backwards, we get that 
$f_{I_{j}}(K) \cap K \neq \emptyset$ for $j<k$. Next, we need to show 
that $I_k$
extends to a sequence in $A$. Since ${\displaystyle f_{I_k}=
\bigcup_j f_{I_k}(f_j(K))}$, there is some $i_{k+1}$ such that 
$f_{i_1}\cdots f_{i_k}f_{i_{k+1}}(K)\cap \partial K \neq \emptyset$. 
Continuing
by induction gives the necessary extension of $I_k$. \qed
\enddemo

Proposition 3.3 leads to a strategy for computing $\dim_H \partial K$.
Suppose for convenience that the $c_i$ all have the same value $c$, and
suppose that we can determine 
${\displaystyle \alpha = \lim_{k \rightarrow \infty} |A_k|^{1/k}}$, where
$|S|$ denotes the cardinality of the set $S$. Then 
$$\tau (s) = \lim_{k \rightarrow \infty} \left( \sum_{I_k \in A_k}
c_{I_k}^s \right)^{1/k}=\lim_{k \rightarrow \infty} 
\left(|A_k|c^{ks}\right)^{1/k}=\alpha c^s.$$ Thus $\dim_H \partial K$ is
the solution to $\alpha c^s=1$, or 
$$\dim_H \partial K = -\frac{\ln(\alpha)}{\ln(c)}.$$

\heading
\S4. The L\'evy Dragon
\endheading

For the rest of this paper, $K$ will denote the L\'evy Dragon, which is the 
attractor of the IFS $F=\{f_1,f_2\}$, where $f_1$ and $f_2$ 
are the similitudes of Euclidean space
$\R^2$ given by
$$\align
 f_1(x,y)&= (\frac{x-y}{2},\frac{x+y}{2}), \\
 f_2(x,y)&= (\frac{x+y+1}{2},\frac{y-x+1}{2}).
\endalign$$ 

L\'evy \cite{L1} studied $K$ extensively, and showed, among other things, that
$K$ tiles the plane in the sense that the plane can be written as the union
of congruent copies of $K$ that meet only in their boundaries. It is first
important that we establish that $K$ satisfies the open set condition so that 
we can use Theorem 2.2.
The following Lemma can be found in \cite{B} or \cite{K}.

\proclaim{4.1 Lemma}  Suppose that $\{ f_1, \dots , f_n \}$ 
is an IFS on $\R^l$ with
$K$ the invariant set.  Suppose that $\sum_{i=1}^n c_i^n = 1$.  Then if
$K$ has nonempty interior in $\R^l$, then $K$ satisfies the open sets
condition with the int $K$ being the open set.
\endproclaim
\demo{Proof}
Let $U = $ int $K$.  Then since each $f_i$ is an open mapping in $R^l$,
it must be the case that $f_i($ int $K) \subset$ int $K$ for each $1 \leq
i \leq l$.  By Lemma 1.1 of \cite{K} it must be the case that
${\Cal H}^l (f_i(K) \cap f_j(K)) = 0$ for all $i \neq j$, where
${\Cal H}^{\alpha}(A)$ denotes the Hausdorff $\alpha$-measure of $A$. 
Since
Lebesgue $l$-dimensional measure is proportional to ${\Cal H}^l$, this
implies that $\lambda (f_i($int $K) \cap (f_j($int $K)) = 0$.  Of course,
this is only possible if $f_i($int $K) \cap (f_j($int $K) = \emptyset$.
\qed
\enddemo

For our study,
we will need some more notation. Let $S$ be the unit square with vertices
$\{(0,0),(0,1),(1,1),(1,0)\}$ and let $T_0 \subset S$ be the triangle with
vertices $\{(0,0),(1/2,1/2),(1,0)\}$. We will view $K$ as 
${\displaystyle \lim_{k \rightarrow \infty} F^k(T_0)}$. 

\epsfysize=2in
\centerline{\epsffile{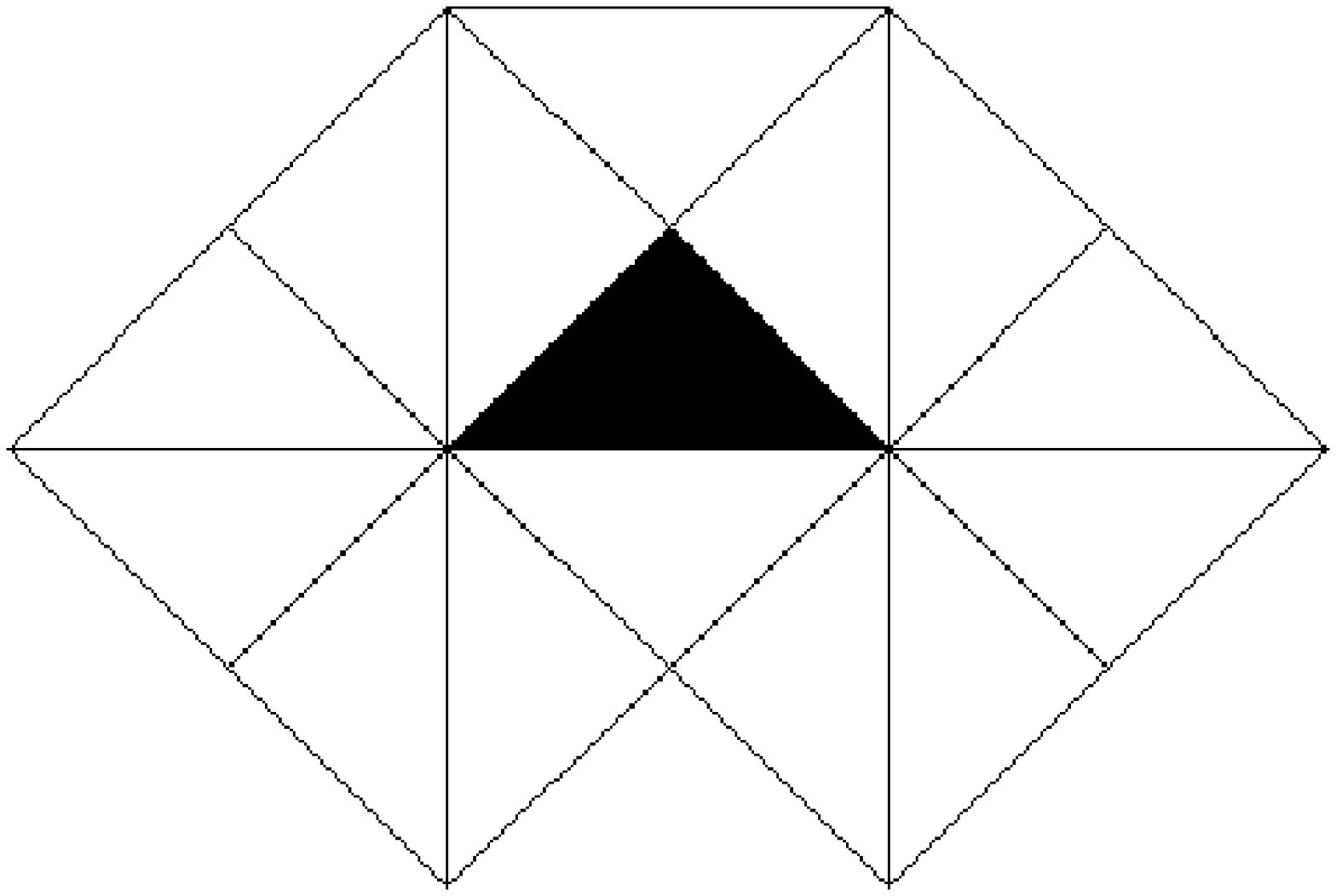}}
\smallskip
\centerline{Figure 2. $T_0$ And Its Neighborhood }
\medskip

$T_0$ and its neighborhood in ${\Cal T}_0$ are shown in Figure 2.
 $S$ is the union of the four triangles in $S$ with
hypotenuse an edge of $S$ and third vertex at $(1/2,1/2)$. We denote the 
triangulation of $\R^2$ consisting of these triangles and their integral
translates by ${\Cal T}_0$. If $T$ is a right triangle in the plane, we
refer to the vertex opposite the hypotenuse as the {\it top vertex} of $T$,
and the other vertices as {\it left} and {\it right} so that we encounter the
vertices in the order left, top, and right as we traverse the boundary 
in clockwise order. The {\it left(right) edge} of $T$ 
is the edge determined by the
left(right) and top vertices of $T$. The triangulation ${\Cal T}_0$ has a
subdivision ${\Cal T}_1$ obtained by subdividing each triangle $T$ into the
two triangles determined by the left(right) edge of $T$ and the midpoint of
the hypotenuse of $T$. Continuing in the obvious way gives a sequence 
${\Cal T}_0 \succ {\Cal T}_1 \succ \cdots \succ {\Cal T}_k\succ \cdots$
of subdivisions of $\R^2$ into right isosceles triangles. Each triangle in
${\Cal T}_k$ has diameter $(\frac{1}{\sqrt2})^{k}$. 
For each $T \in {\Cal T}_k$, let $T^L$ be the triangle in ${\Cal T}_{k+1}$
whose hypotenuse is the left edge of $T$ and which has a vertex outside of $T$.
Define $T^R$ similarly. One can easily see that $F(T_0)=T_0^L \bigcup T_0^R$.
In fact, the following Proposition is immediate.

\proclaim{4.2 Proposition} For each $k>0$, $F^k(T_0)$ is a union of $2^k$ 
triangles in ${\Cal T}_k$. $F^{k+1}(T_0)$ is obtained from $F^k(T_0)$ by
replacing each $T$ in $F^k(T_0)$ by $T^L$ and $T^R$.  \qed
\endproclaim

We would like to use the formula at the end of \S3 to compute the 
Hausdorff dimension of $\partial K$. The problem is that $K$ and 
$\partial K$ are very complicated objects, and we can only see finite
approximations to them, so it is not clear how to compute $A_k$. What we
can hope to compute is $B_k=\{I_k \in \Omega_k|f_{i_k}(T_0) \cap
\partial F^k(T_0) \neq \emptyset\}$. We now attack that problem together
with the question of relating $A_k$ to $B_k$.

For $T \in {\Cal T}_k$, let $N(T)$ be the star of $T$, that is, 
$N(T)=\{T' \in {\Cal T}_k|T \cap T' \neq \emptyset\}$. Each $N(T)$ consists
of 15 triangles. For technical reasons to be seen later, we order these 
triangles as $N(T)[1],N(T)[2],\cdots ,N(T)[15]$ as follows: $T=N(T)[1]$;
then $N(T)[2]-N(T)[8]$ are the triangles incident with the left vertex of $T$ in
clockwise order, so that $N(T)[8]$ has the left edge of $T$ for its right edge;
$N(T)[9]$ is incident with the top vertex of $T$ and and has the left edge of
$N(T)[8]$ for its right edge;$N(T)[10]$ has the right edge of $T$ for its
left edge; and $N(T)[11]-N(T)[15]$ are the remaining triangles incident
with the right vertex of $T$ ordered clockwise. 
This ordering is pictured in Figure 3.
The {\it neighborhood
type} $\rho(T)$ is the 15-tuple $\rho(T)=\{\rho_1,\rho_2,\cdots,\rho_{15}\}$,
where $$
\rho _i=\cases 1,&\text{if $N(T)[i] \in F^k(T_0);$}\\
              0,&\text{otherwise.}\endcases$$
We say that $T \in {\Cal T}_k$ is of {\it type\/} $\rho(T)$ and that $T$ 
is {\it covered} if $\rho(T)$ is the all ones
vector. $T$ being covered means that $T$ and all of the triangles in 
${\Cal T}_k$
that intersect $T$ are contained in $F^k(T_0)$. To avoid even further notation,
we will use $N(T)$ to denote both the collection of triangles as defined above 
as well as the set which is the union of those triangles (the {\it carrier}
of $N(T)$).

\epsfysize=2in
\centerline{\epsffile{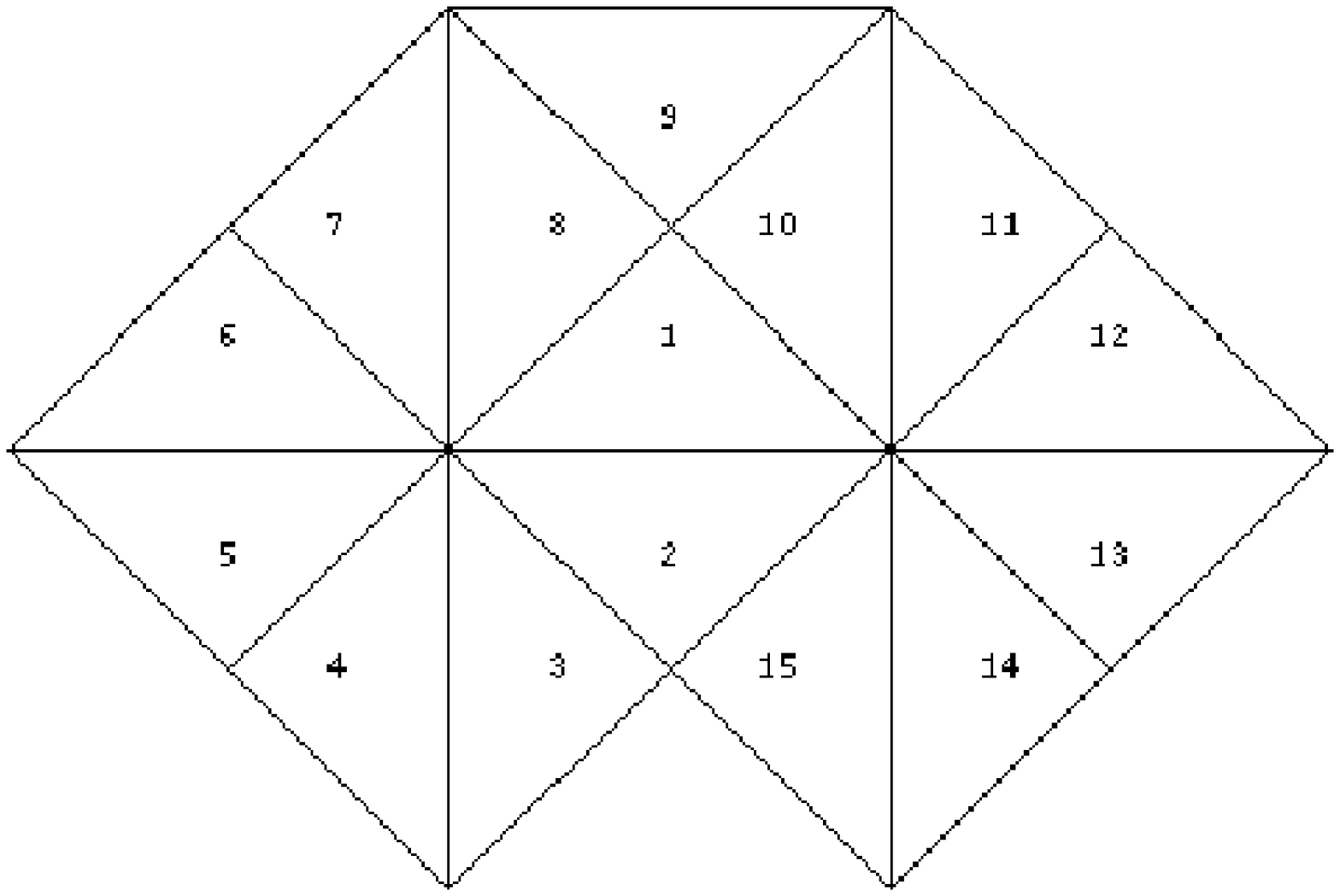}}
\smallskip
\centerline{Figure 3. The Neighborhood Ordering}
\medskip

One of the difficulties in dealing with $K$ is that it is not easy to locate
points in the interior of $K$ by looking at $F^k(T_0)$. One needs to find
objects of positive area which persist from iteration to the next, but as
soon as a triangle $T$ appears at one level, its interior is discarded at the 
next level and is replaced by two triangles exterior to $T$. It is not until
the iterations have become sufficiently intertwined to create covered triangles
that the discarded mass is filled back in by adjacent triangles. Through 
computer experimentation we discovered that covered
triangles do not exist until the 14th iteration, and only 8 of the $2^{14}$
triangles that comprise $F^{14}(T_0)$ are covered!
However, if a triangle 
$T \in {\Cal T}_k$
is covered, then $T$ is the union of two triangles in ${\Cal T}_{k+1}$ which
are also covered. It follows that int $T \subset F^{k+m}(T_0)$ for all
$m>0$, so that int $T \subset$ int $K$. We record this and
a few other elementary observations about neighborhoods in the next
 Proposition.

\epsfysize=2in
\centerline{\epsffile{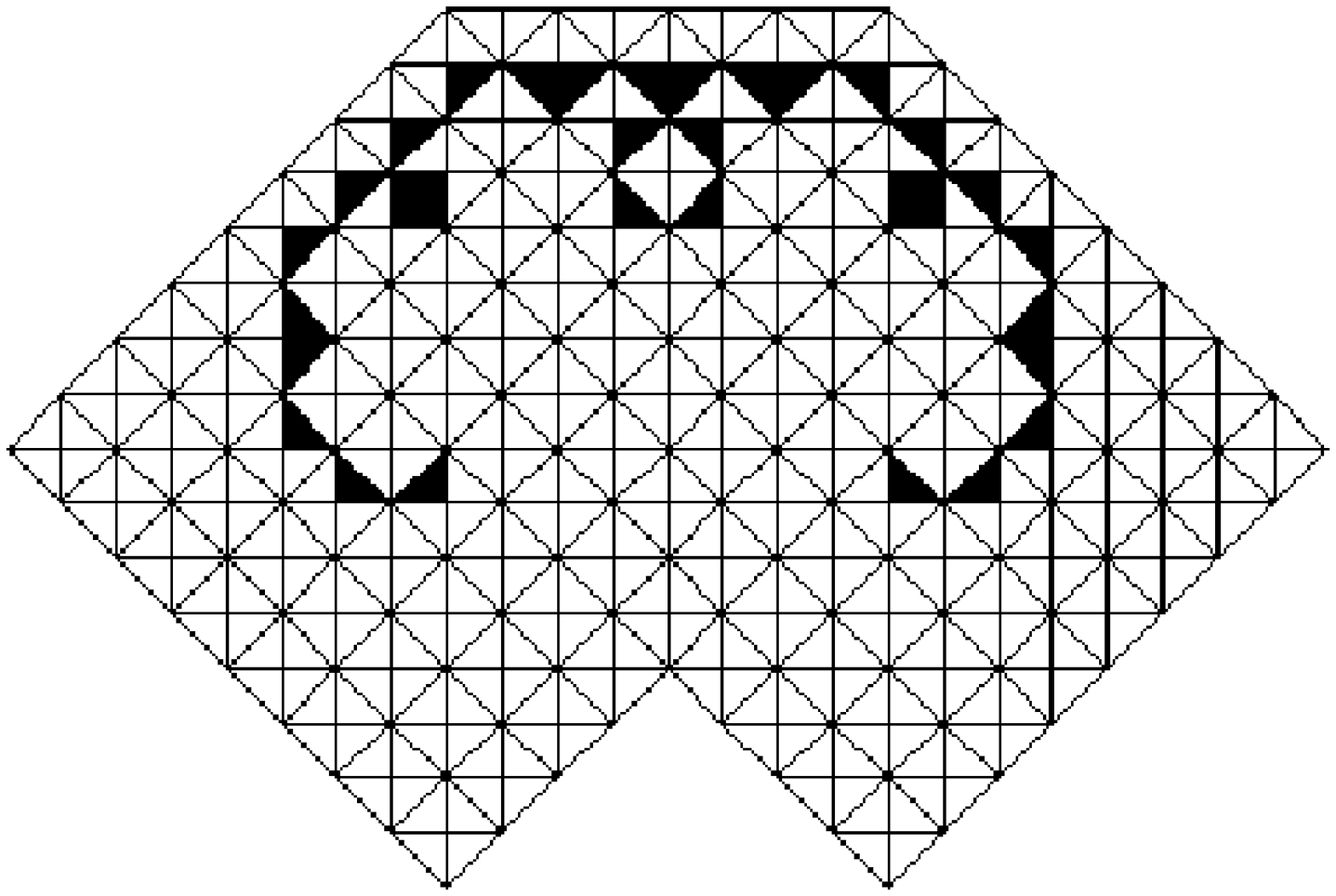}}
\smallskip
\centerline{Figure 4. $N_5$ And $F^5{T_0}$}
\medskip

\proclaim{4.3 Proposition}
\roster
\item Let $N_0=N(T_0)$. Then $F(N_0) \subset N_0$, and therefore
$K \subset N_0$.
\item If $T \in {\Cal T}_k$ and $T' \subset T$ for some $T'\in {\Cal T}_{k+1}$,
then $N(T') \subset N(T)$.
\item If $T \in {\Cal T}_k$ is covered and $T' \subset T$ for some 
$T'\in {\Cal T}_{k+1}$, then $T'$ is covered.
\item If $T \in {\Cal T}_k$ is covered, then ${\text int }T 
\subset {\text int }K$. 
\item If $v$ is a vertex of any triangle $f_{I_k}(T_0) \subset F^k(T_0)$, it
 is also a vertex of some triangle $f_{I_{k+1}}(T_0) \subset F^{k+1}(T_0)$,
 and therefore $v \in K$. \qed
\endroster
\endproclaim

As we mentioned above, L\'evy [L1] observed that the 
Dragon $K$ tiles the plane, and
the tiles occur in the following way. For each triangle $T_a \in {\Cal T}_0$
there is an $IFS$ $F_a$ whose similarities are related to $T_a$ in the same
way that $f_1$ and $f_2$ are related to $T_0$. At each stage, the 
$F_a^k(T_a)$ tile the plane, and the limits $\{K_a\}$ give a tiling by copies
of the Dragon. 

Now recall our definitions $B_k=\{I_k \in \Omega_k|f_{I_k}(T_0) \cap
  \partial F^k(T_0) \neq \emptyset\}$ and
$A_k=\{I_k \in \Omega_k|f_{I_k}(K) \cap \partial K \neq \emptyset\}$.

\proclaim{4.4 Proposition} $B_k \subset A_k$. \endproclaim
\demo{Proof} 
Note that $B_k$ is the set of $I_k$ such that $f_{I_k}(T_0)$ is not covered.
Suppose that $I_k \in B_k$, and let $T=f_{I_k}(T_0)$. Then for some $T' \in
N(T)$, $T' \nsubseteq F^k(T_0)$ and thus $T' \subset F_a^k(T_a)$ for some
other tile at the $k$th level. Let $v$ be a vertex common to $T$ and $T'$.
By 4.3 (5) applied to both $F$ and $F_a$, $v \in K \cap K_a$, so
$v \in \partial K$. Since $v = f_{I_k}(w)$ for some vertex $w$ of $T_0$,
$v \in f_{I_k}(K)$, so $f_{I_k}(K) \cap \partial K \neq \emptyset$, and
$I_k \in A_k$. 
\enddemo

It follows that $|B_k| \leq |A_k|$. By using $B_n$ in two ways, we can compute
$\dim_H \partial K$. Let ${\displaystyle \beta = \lim_{k \rightarrow \infty}
|B_k|^{1/k}}$.

\proclaim{4.5 Theorem} $$\dim_H \partial K = \frac{\ln(\beta)}{\ln(\sqrt2)}.$$

Moreover, $0 < H^\alpha(\partial K) < +\infty$ where $H^\alpha(A)$ is the 
Hausdorff $\alpha$-dimensional measure and $\alpha = \frac{\ln(\beta)}{\ln(\sqrt2)}$.
\endproclaim

\demo{Proof}
Since $\beta \leq \alpha$, we get $$\frac{\ln(\beta)}{\ln(\sqrt2)} \leq
\dim_H \partial K$$ from the formula at the end of \S3. On the 
other hand, we can also use $\beta$ to get an upper bound for the box counting
dimension of $\partial K$, which is equal to the Hausdorff
dimension by Theorem 2.2. Let $N_k$ denote the  
the set of triangles $T$ in the subdivision
of $N_0$ induced by ${\Cal T}_k$. For such a $T$, it is clear that if
$\rho(T)$ is the all zeros vector, then $\text{int }T \cap K = \emptyset$ 
and if $T$ is covered, $\text{int }T \cap \partial K = \emptyset$. Thus,
$\partial K \subset \bigcup\{T \in N_k|\rho(T) \text{ is not constant}\}$. By
expanding to neighborhoods, we get 
$\partial K \subset \bigcup\{N(T)| N(T) \text{ contains some } f_{I_k}(T_0),
I_k \in B_k\}$. There are at most $15|B_k|$ such $N(T)$, so for each $k$,
$\partial K$ is covered by at most $15|B_k|$ sets of diameter 
$\leq 3(\sqrt 2)^{-k}$. Therefore, the box counting dimension of $\partial K$
is less than or equal to 
$$\lim_{k \rightarrow \infty} 
-\frac{\ln(15|B_k|)}{\ln(3(\sqrt 2)^{-k})} 
= \frac{\ln(\beta)}{\ln(\sqrt2)}.$$

The referee pointed out the second part of the theorem.  It is one of the
advantages of the approach of Strichartz and Wang \cite{SW} using the graph
directed construction of Mauldin and Williams \cite{MW} that one can claim
that $0 < H^\alpha(\partial K) < +\infty$.  It also follows from 
\cite{MW} and the fact in this paper that the matrix $M$ used in arriving
at the asymptotic estimate of $|B_k|$ is irreducible. \qed
\enddemo

\heading
\S5. Computational issues
\endheading

In the light of Theorem 4.5, our goal is to understand the asymptotic
behavior of the number
$|B_k|$ of triangles in $F^k(T_0)$ which are not covered. To this end, we 
need to study the neighborhood structures introduced in \S4. We will index
the binary vectors of length 15 by integers between 0 and $2^{15}-1 = 32767$
by the correspondence $(x_1,x_2,\cdots,x_{15}) \Leftrightarrow 
\sum_{i=1}^{15}x_i2^{i-1}$. Thus we could indicate that $T$ is covered by
writing $\rho(T)=32767$. Recall that each $T \in {\Cal T}_k$ is the union of 
two triangles $T_1$ and $T_2$ in ${\Cal T}_{k+1}$, where we let $T_1$ be the
triangle that contains the left vertex of $T$, and $T_2$ contains the right
vertex. If we know the neighborhood type for $T$, it is simple (but
tedious) to write down the neighborhood types for $T_1$ and $T_2$. We
have

\proclaim{5.1 Proposition} For $T \in {\Cal T}_k$, if
$$\rho(T)=(x1, x2, x3, x4, x5, x6, x7, x8, x9, x10, x11, x12, x13, x14, x15),$$
then 
$$\rho(T_1)=(x8, x1, x9, x8, x10, x9, x1, x10, x15, x3, x2, x5, x4, x7, x6),$$
and
$$\rho(T_2)=(x10,x1,x12,x11,x14, x13, x2, x15, x3, x8, x1, x9, x8, x10, x9).
\qed$$
\endproclaim

Proposition 5.1 gives a way to compute $|B_k|$ in principle. Let $V(k)$ be
the $2^{15}$ long vector such that $V(k)_i = \text{ number of } T \in N_k 
\text{ such that } \rho(T)=i$. Then let M be the $2^{15} \times 2^{15}$
matrix with $M_{i,j}= t$, if a triangle of type $i$ gives
rise to $t$ subtriangles with type $j$. For example, 
$M_{0,0}=2,M_{32767,32767}=2, M_{0,j}=0,i>0$,and $M_{32767,j}=0,j<32767$.
$M$ can be constructed from the information in 5.1. Then $M$ has the
property that $V(k+1)=V(k)\cdot M$, so if $J$ is the column vector with
$J_i = 1$ if $i$ is odd and less than 32767, 0 otherwise, then
$|B_k| = V(0)\cdot M^k \cdot J$ and we could study the 
asymptotic behavior of $|B_k|$
via the eigenvalues of $M$. Fortunately, we do not have to go to such 
extremes. A little thought shows that not every one of the $2^{15}$
neighborhood types can occur.  Our computations determined that, in fact,
relatively few types actually do occur.

The only triangle in $F^0(T_0)$ is $T_0$ itself, so the neighborhood 
types for the fifteen triangles in $N_0$ are encoded by $1,2,2^2,
\cdots,2^{14}$. For example, $T_0$ is triangle 1 in its neighborhood,
so $\rho(T_0)=1$. For a more revealing example, consider the triangle $T$ with
vertices $\{(0,0),(0,1),(1/2,1/2)\}$. It is number 8 in the neighborhood
of $T_0$, but $T_0$ is number 10 in $N(T)$. Therefore, $\rho(T)=2^{(10-1)}
=512$. If we let $S_k = \{\rho(T)|T \in N_k\}$, then, $S_0=\{2^k,0 \leq
k \leq 14\}$. By applying the formulas above for $\rho(T_1)$ and $\rho(T_2)$,
we can compute $S_1,S_2,$ and so on.
Again, we determined by computation that $S_{19} = S_{20}$,
so that $S_{k+19}=S_{19}$ for all $k>0$. 
Let $S_{\infty}$ be this common ``stable''
set of neighborhood types, $S_{\infty}=\{s_1,s_2,\cdots,s_{752}\}$.
(We will discuss 
the order on $S_\infty$ below.)
The cardinality of $S_{\infty}$ is 752, so
we can perform the analysis indicated above with state vector and transition
matrix of manageable size. For $k \geq 19$, let $V(k)$ be the vector 
of length 752 whose $i$th entry is the number of triangles in $N_k$ which have
type $s_i$.  Let $M$ be the 
$752 \times 752$ matrix with $M_{i,j}= t$, if a triangle of type $s_i$ gives
rise to $t$ subtriangles with type $s_j$. 
Then as above,$V(k+1)=V(k)\cdot M$, and
if $J$ is the column vector whose $i$th entry is 1 if $s_i$ is coded by
an odd number $\neq 32767$, and 0 otherwise, then
 $$|B_k| = V(k)\cdot M \cdot J,$$ so that
$$|B_{19+k}| = V(19)\cdot M^k \cdot J.$$ 

The types in $S_{\infty}$ have a nice structure. The types 0 and 32767 are
absorbing in that each gives rise to two triangles of the same type at the
next level. Let $S_a=\{0,32767\}$; let
$$S_t=\{4,8,9,16,32,64,66,128,512,1024,1026,2048,4096,8192,8193,16384\};$$
and let $S_e = S_{\infty} \setminus (S_a \bigcup S_t)$. It can be shown
by computation that if we order the types in $S_{\infty}$ so that 
$s_i < s_j < s_k$ if $s_i \in S_t, s_j \in S_e,s_k \in S_a$, then $M$ has the
block form 
$$ M = \left( \matrix
P & Q & R\\
0 & C & L\\ 
0 & 0 & I
\endmatrix \right),$$
where $C$ is a $734\times734$ matrix and $I$ is the $2\times2$ identity matrix.
Write 
$V(k)=(V_t(k),V_e(k),V_a(k))$, 
and let 
$$J=\left( \matrix
J_t\\
J_e\\
J_a\endmatrix \right),$$
where the blocking is consistent with that of
$M$.
 Then we have
$$ |B_{19+k}| =(V_t(k),V_e(k))\cdot \left( \matrix
P & Q\\
0 & C \endmatrix \right)^k \cdot \left( \matrix
J_t\\
J_e\endmatrix \right).$$

The matrix $P$ can be seen to be
a permutation matrix, and $C$ is primitive in the sense
of Perron-Frobenius theory \cite{BP}, {\it i.e.\/}, all entries of some
power of $C$ are positive. In fact, a calculation shows that $C^{25}\gg 0$.
It follows \cite{BP, Theorem 2.4.1} that $C$ has an eigenvalue $\lambda$
equal to its spectral radius, and that ${\displaystyle 
\lim_{k \rightarrow \infty}\frac{1}{\lambda^k}C^k = D}$, where $D$ is a 
positive matrix. In our case, $\lambda>1$, so 
$$\lim_{k \rightarrow \infty}\frac{1}{\lambda^k}
\left( \matrix
P & Q\\
0 & C \endmatrix \right)^k =
\lim_{k \rightarrow \infty}\left( \matrix
\frac{1}{\lambda^k}P^k & \frac{1}{\lambda^k}P^{k-1}Q\\
0 & \frac{1}{\lambda^k}C^k \endmatrix \right) =
\left( \matrix
0 & 0\\
0 & D \endmatrix \right).$$ 
Therefore,
$$\lim_{k \rightarrow \infty} \frac{|B_{k+19}|}{\lambda^{k+19}}=
\frac{1}{\lambda^{19}}V_e(19) \cdot D \cdot J_e = q > 0.$$
Then we have $$\beta = \lim_{m \rightarrow \infty}|B_m|^{1/m}=
\lim_{m \rightarrow \infty}\lambda q^{1/m} = \lambda.$$

By \S 4, $$\dim_H \partial K = \frac{\ln(\lambda)}{\ln(\sqrt 2)}.$$

Floating point computations using the power method give the estimate
$\lambda \approx
1.954776399$.  That approximation for $\lambda$ gives $\dim_H
\partial K \approx 1.934007183$.
According to \cite{BP, Theorem 2.2.35}, the spectral radius 
$r$ of an irreducible
matrix $L$ satisfies $u \leq r \leq U$, where $u$ and $U$ are the minimum
and maximum row sums of $L$, so it follows that our $\lambda$ satisfies
$u_k^{1/k} \leq \lambda \leq U_k^{1/k}$ for each $k$, 
where $u_k$ and $U_k$ are the minimum and maximum row sums of $C^k$. We 
computed $u_{30}=67936360$ and $U_{30}=727212953$
using integer arithmetic, which gives the rigorous
estimate $1.824190 < \lambda < 1.974189$. Thus $\lambda$ is provably greater
than $\sqrt 2$, and  the Hausdorff dimension of $\partial K$ is greater than
1.  Thus, the boundary of the Dragon has Hausdorff dimension greater than 1,
and Edgar's question is answered in the affirmative.

\enddocument